\documentclass[12pt]{article}
\usepackage{algorithm,amsfonts,amsmath,amssymb,bbm,fullpage,graphicx,varwidth}
\usepackage{mathtools,psfrag,verbatim}
\usepackage{algpseudocode}
\usepackage[small,bf]{caption}
\usepackage{subcaption, color}
\usepackage{url}

\newcommand{\reals}{{\mbox{\bf R}}}

\newcommand{\Tr}{\mathop{\bf Tr}}

\newcommand{\ie}{\emph{i.e.}}

\newcommand{\Split}{\mathop{\mathrm{Split}}}

\newcommand{\Kmax}{K^\mathrm{max}}

\newcommand{\BEAS}{\begin{eqnarray*}}
\newcommand{\EEAS}{\end{eqnarray*}}
\newcommand{\BEQ}{\begin{equation}}
\newcommand{\EEQ}{\end{equation}}
\newcommand{\BIT}{\begin{itemize}}
\newcommand{\EIT}{\end{itemize}}

\title{Greedy Gaussian Segmentation\\ of Multivariate Time Series}
\author{David Hallac \and Peter Nystrup \and Stephen Boyd} 
\date{\small{April 2018}}

\begin{document}
\maketitle

\makeatletter
\def\BState{\State\hskip-\ALG@thistlm}
\makeatother

\algnewcommand\algorithmicinput{\textbf{Input:}}
\algnewcommand\Input{\item[\algorithmicinput]}

\begin{abstract}
We consider the problem of breaking a multivariate (vector) time series into
segments over which the data is well explained as independent samples 
from a Gaussian distribution.
We formulate this as a covariance-regularized maximum likelihood problem,
which can be reduced to a combinatorial optimization problem of
searching over the possible breakpoints, or segment boundaries.
This problem can be solved using dynamic programming, with complexity
that grows with the square of the time series length.
We propose a heuristic method that approximately solves the problem 
in linear time with respect to this length,
and always yields a locally optimal choice, in the sense that no change of
any one breakpoint improves the objective.
Our method, which we call \emph{greedy Gaussian segmentation} (GGS),
easily scales to problems with
vectors of dimension over 1000 and time series of arbitrary length.
We discuss methods that can be used to validate such a model
using data, and also to automatically choose appropriate values of the two 
hyperparameters in the method.
Finally, we illustrate our GGS approach on financial time series and Wikipedia text data.
\end{abstract}


\section{Introduction} 

Many applications, including weather measurements \cite{xu2002survey}, 
car sensors \cite{hallac2016driver}, and financial returns 
\cite{nystrup2017dynamic}, 
contain long sequences of multivariate time series data. 
With data sets such as these, there are many benefits to 
partitioning the time series into segments, where each segment
can be explained by as simple a model as possible. 
Partitioning can be used for denoising \cite{abonyi2005modified}, 
anomaly detection \cite{rajagopalan2006symbolic}, 
regime-change identification \cite{nystrup2016detecting}, and more.
Breaking a large data set down into smaller, simpler components is also
a key aspect of many 
unsupervised learning algorithms \cite[Chapter 14]{hastie2009elements}.

In this paper, we analyze the time series partitioning problem
by formulating it as a covariance-regularized 
likelihood maximization problem, where the data in each segment can 
be explained as independent samples from a multivariate Gaussian 
distribution. We propose an efficient heuristic, which we call
the \emph{greedy Gaussian segmentation} (GGS) algorithm, that
approximately finds the optimal breakpoints using a greedy
homotopy approach based on the number of segments
\cite{zangwill1981pathways}.
The memory usage of the algorithm is a modest multiple of the 
memory used to represent the original data,
and the time complexity is linear in the number of 
observations, with significant opportunities for exploiting parallelism.
Our method is able to scale to arbitrarily long time series
and multivariate vectors of dimension over 1000. We also discuss
several extensions of this approach, including a streaming algorithm
for real-time partitioning, as well as a method of validating the model
and selecting optimal values of the hyperparameters. 
Last, we implement the GGS algorithm in a Python software package \texttt{GGS}, available
online 
at \url{https://github.com/cvxgrp/GGS}, and apply it to various financial time 
series and Wikipedia text data to illustrate our method's accuracy, scalability, 
and interpretability.

\subsection{Related work}

This work relates to recent advancements in both optimization and 
time series segmentation.
Many variants of our problem
have been studied in several contexts, including Bayesian change-point detection
\cite{booth1982bayesian, lee1998bayesian, son2005bayesian,
cheon2010multiple, bauwens2012marginal},
change-point detection based on hypothesis testing
\cite{crosier1988multivariate, venter1996finding, de2006detecting,
galeano2014multiple, li2015nonparametric},
mixture models \cite{verbeek2003efficient, abonyi2005modified, 
picard2011joint, same2011model},
hidden Markov models and the Viterbi algorithm 
\cite{ryden1998stylized, ge2001segmental, bulla2011hidden,
hu2015using, nystrup2017long}, and
convex segmentation \cite{katz2014outlier}, all trying to find 
breakpoints in time series data sets.

The different methods make different assumptions about the
data (see \cite{esling2012time} for a comprehensive survey). GGS 
assumes that, in each segment,
the mean and covariance are constant and unrelated to the means
and covariances in all other segments. 
This differs from ergodic hidden Markov models
\cite{ryden1998stylized, ge2001segmental, bulla2011hidden,
hu2015using, nystrup2017long},
which implicitly assume that the underlying segments will
repeat themselves, with some structure to when the transitions are likely 
to occur. In a left-to-right hidden Markov model \cite{bakis1976continuous, cappe2005inference},
though, additional constraints are imposed to ensure non-repeatability of segments,
similar to GGS. Alternatively, trend filtering problems 
\cite{kim2009L1trend} assume that neighboring segments have similar
statistical parameters; when a transition occurs, the new parameters
are not too far from the previous ones.
Other models have tried to solve the problem of change-point detection 
when the number of 
breakpoints is unknown \cite{basseville1993detection, chouakria2003compression},
including in streaming settings \cite{guralnik1999event, gustafsson2000adaptive}.

GGS uses a straightforward
approach based on the maximum likelihood of the data (we address how to
incorporate many of these alternative assumptions in \S\ref{s-variations}).
In real world contexts, deciding on which approach to use depends entirely
on the underlying structure of the data; a reasonable choice of 
method can be determined via cross-validation of the various models.
Our work is novel in that it allows for an extremely scalable 
greedy algorithm to detect breakpoints in multivariate time series. 
That is, GGS is able to
solve much larger problems than many of these other methods, 
both in terms of vector dimension and the length of the time series.
Additionally, its robustness allows GGS to be used as a black-box method
which can automatically determine an appropriate number of 
breakpoints, as well as the model parameters within each segment,
using cross-validation.

Our greedy algorithm is based on a top-down approach to segmentation	
\cite{douglas1973algorithms}, though there has
also been related work using bottom-up methods \cite{keogh2004segmenting}.
While our algorithm does achieve a locally optimal solution, we note that it is
possible to solve for the global optimum using dynamic programming
\cite{bellman1961approximation, fragkou2004dynamic, kehagias2006dynamic}.
However, these globally optimal approaches have complexities that
grow with the square of the time series length, whereas our heuristic method
scales linearly with the time series length.
Our model approximates
$\ell_1$/$\ell_2$ trend filtering problems \cite{kim2009L1trend, 
wahlberg2011L1mean, wahlberg2012ADMM},
which use a penalty based on the fused group lasso
\cite{tibshirani2005sparsity, bleakley2011group}, to 
couple together the model parameters at adjacent times.
However, these models are unable to scale up
to the sizes we are aiming for, so we develop a fast heuristic, similar to an 
$\ell_0$ penalty \cite{candes2008enhancing}, where each breakpoint splits 
the time series into two independent problems.
To ensure robustness, we rely on covariance-regularized 
regression to avoid errors when there are more
dimensions than samples in a segment \cite{witten2009covariance}.

%Our greedy algorithm is based on a top-down approach to segmentation	
%\cite{douglas1973algorithms}, though there has
%also been related work using bottom-up methods \cite{keogh2004segmenting}.
%While some combinatorial optimization problems have an efficient globally 
%optimal solution method based on dynamic programming 
%\cite{bellman1957dynamic, bertsekas1995dynamic}, the problem that GGS attempts 
%to solve cannot be solved by such a method, to the best of our knowledge.
%Our model approximates
%$\ell_1$/$\ell_2$ trend filtering problems \cite{kim2009L1trend, 
%wahlberg2011L1mean, wahlberg2012ADMM},
%which typically use a penalty based on the fused group lasso
%\cite{tibshirani2005sparsity, bleakley2011group}, to 
%couple together the model parameters at adjacent times.
%However, these models are unable to scale up
%to the sizes we are aiming for, so we develop a fast heuristic, similar to an 
%$\ell_0$ penalty \cite{candes2008enhancing}, where each breakpoint splits 
%the time series into two independent problems.
%To ensure robustness, we rely on covariance-regularized 
%regression to avoid errors when there are more
%dimensions than samples in a segment \cite{witten2009covariance}.

\subsection{Outline}

The rest of this paper is structured as follows. In \S\ref{s-problem}, we formally
define our optimization problem. In \S\ref{s-alg}, we explain the GGS
algorithm for approximately solving the problem in a scalable way.
In \S\ref{s-validation}, we describe a validation process 
for choosing the two hyperparameters in our model.
We then examine in \S\ref{s-variations} several extensions of this approach which
allow us to apply our algorithm to new types of problems. 
Finally, we apply GGS to several real-world financial and Wikipedia
data sets, as well as a synthetic example, in \S\ref{s-examples}.
%Finally, we conclude and list potential 
%directions for future work in \S\ref{s-conclusion}.

\section{Problem setup}\label{s-problem}
\subsection{Segmented Gaussian model}
We consider a given time series $x_1, \ldots, x_T \in \reals^n$.
(The times $t=1,\ldots, T$ need not be uniformly spaced in real time;
all that matters in our model and method is that they are ordered.)
We will assume that the $x_t$'s are independent samples with $x_t
\sim\mathcal{N}(\mu_t, \Sigma_t)$,
where the mean $\mu_t$ and covariance $\Sigma_t$ only change at $K\ll T$
breakpoints $b_1, \ldots, b_K$.
These breakpoints divide the given $T$ samples into $K+1$ 
segments;
in each segment, the $x_t$'s are generated from the same multivariate Gaussian distribution.
Our goal is to determine $K$, the break points 
$b_1,\ldots, b_K$, and the means and covariances 
\[
\mu^{(1)}, \ldots, \mu^{(K+1)}, \qquad
\Sigma^{(1)}, \ldots, \Sigma^{(K+1)}
\]
in the $K+1$ segments between the breakpoints,
from the given data $x_1, \ldots, x_T$.

Introducing breakpoints $b_0$ and $b_{K+1}$,
the breakpoints must satisfy
\[
1 = b_0 <b_1< \cdots < b_K < b_{K+1} = T+1,
\]
and the means and covariances are given by
\[
(\mu_t,\Sigma_t) = (\mu^{(i)},\Sigma^{(i)}), \quad
b_{i-1} \leq t < b_i, \quad i=1, \ldots, K.
\]
(The subscript $t$ denotes time $t$; the superscript $(i)$ and subscript
on $b$ denotes segment $i$.)
We refer to this parametrized distribution of $x_1, \ldots, x_T$
as the \emph{segmented Gaussian model} (SGM).
The log-likelihood of the data $x_1, \ldots, x_T$ under this model is
given by
\BEAS
\ell(b, \mu, \Sigma) 
&=&
\sum_{t = 1}^T \left( 
-\frac{1}{2} (x_t - \mu_t)^T \Sigma_t^{-1} (x_t - \mu_t) - 
\frac{1}{2} \log\det \Sigma_t - \frac{n}{2}\log(2\pi) \right) \\
&=& 
\sum_{i=1}^{K+1}
\sum_{t = b_{i-1}}^{b_i-1} \left( 
-\frac{1}{2} (x_t - \mu^{(i)})^T (\Sigma^{(i)})^{-1} (x_t - \mu^{(i)}) - 
\frac{1}{2} \log\det \Sigma^{(i)} - \frac{n}{2}\log(2\pi) \right)\\
&=& 
\sum_{i=1}^{K+1} \ell^{(i)} (b_{i-1}, b_i, \mu^{(i)}, \Sigma^{(i)}),
\EEAS
where
\BEAS
\ell^{(i)} (b_{i-1}, b_i, \mu^{(i)}, \Sigma^{(i)}) &=&
\sum_{t = b_{i-1}}^{b_i-1} \left( 
-\frac{1}{2} (x_t - \mu^{(i)})^T (\Sigma^{(i)})^{-1} (x_t - \mu^{(i)}) - 
\frac{1}{2} \log\det \Sigma^{(i)} - \frac{n}{2}\log(2\pi) \right)\\
&=&
-\frac{1}{2} \sum_{t = b_{i-1}}^{b_i-1} 
(x_t - \mu^{(i)})^T (\Sigma^{(i)})^{-1} (x_t - \mu^{(i)}) \\ 
&& \quad
-\frac{b_i-b_{i-1}}{2} \left( \log\det \Sigma^{(i)} + n \log(2\pi) \right)
\EEAS
is the contribution from the $i$th segment.
Here we use the notation $b=(b_1, \ldots, b_K)$,
$\mu = (\mu^{(1)}, \ldots, \mu^{(K+1)})$, and
$\Sigma = (\Sigma^{(1)}, \ldots, \Sigma^{(K+1)})$, for the parameters in the SGM.
In all the expressions above we define $\log \det \Sigma$ as 
$-\infty$ if $\Sigma$ is singular, \ie, not positive definite. 
Note that $b_i - b_{i-1}$ is the length of the $i$th segment.

\subsection{Regularized maximum likelihood estimation}\label{s-regularized}
We will choose the model parameters by maximizing the 
covariance-regularized log-likelihood for a given value of $K$, 
the number of breakpoints. We regularize the covariance 
to avoid errors when there are more dimensions than samples in a segment,
a well-known problem in high dimensional settings
\cite{huang2006covariance, bickel2008regularized, witten2009covariance}.
Thus we choose $b,\mu,\Sigma$ to maximize the regularized log-likelihood
\BEQ\label{e-reg-ll}
\phi(b,\mu,\Sigma) = 
\ell(b,\mu,\Sigma) - \lambda \sum_{i=1}^{K+1} \Tr (\Sigma^{(i)})^{-1}
= \sum_{i=1}^{K+1} \left( 
\ell^{(i)}(b_{i-1},b_i,\mu^{(i)},\Sigma^{(i)}) - 
\lambda \Tr (\Sigma^{(i)})^{-1}\right),
\EEQ
where $\lambda \geq 0$ is a regularization parameter, with $K$ fixed.
(We discuss the choice of the hyperparameters 
$\lambda$ and $K$ in \S\ref{s-validation}.)
This is a mixed combinatorial and continuous 
optimization problem since it involves a
search over the $T-1 \choose K$ possible choices of the breakpoints 
$b_1, \ldots, b_K$, as well as the parameters $\mu$ and $\Sigma$.
For $\lambda =0$, this reduces to maximum likelihood estimation, 
but we will assume henceforth that $\lambda>0$.
This implies that we will only consider positive definite (invertible) 
estimated covariance matrices.

If the breakpoints $b$ are fixed, the regularized maximum likelihood 
problem has a simple analytical
solution.  The optimal value of the $i$th segment mean is the empirical 
mean over the segment,
\BEQ \label{e-emp-mean}
\mu^{(i)} = \frac{1}{b_i - b_{i-1}}\sum_{t = b_{i-1}}^{b_i-1} x_t,
\EEQ
and the optimal value of the $i$th segment covariance is
\BEQ\label{e-cov}
\Sigma^{(i)} = S^{(i)} + \frac{\lambda}{b_i - b_{i-1}} I,
\EEQ
where $S^{(i)}$ is the empirical covariance over the segment,
\[
S^{(i)} = \frac{1}{b_i - b_{i-1}}\sum_{t = b_{i-1}}^{b_i-1}
(x_t-\mu^{(i)})(x_t-\mu^{(i)})^T.
\]
Note that the empirical covariance $S^{(i)}$ can be singular, for example 
when $b_i-b_{i-1}<n$, but
for $\lambda >0$ (which we assume), 
$\Sigma^{(i)}$ is always positive definite.
Thus, for any fixed choice of breakpoints $b$, the mean and covariance parameters 
that maximize the regularized log-likelihood \eqref{e-reg-ll}
are given by \eqref{e-emp-mean} and \eqref{e-cov}, respectively.
The optimal value of the covariance \eqref{e-cov} is similar to a
Stein-type shrinkage estimator \cite{ledoit2004well}.

Using these optimal values of the mean and covariance parameters,
the regularized log-likelihood \eqref{e-reg-ll} can be expressed in terms of 
$b$ alone, as
\BEAS
\phi(b) %&=& -\frac{Tn}{2}\log(2\pi) -\frac{1}{2} \sum_{i=1}^{K+1} 
%(b_i-b_{i-1}) 
%\left( \Tr (S^{(i)}+\frac{\lambda}{b_i - b_{i-1}} I)^{-1} S^{(i)}
%+ \log\det (S^{(i)}+\frac{\lambda}{b_i - b_{i-1}} I) \right) \\
&=& C - \frac{1}{2} \sum_{i=1}^{K+1} \left(
(b_i-b_{i-1}) \log\det (S^{(i)}+\frac{\lambda}{b_i - b_{i-1}} I)
- \lambda \Tr (S^{(i)}+\frac{\lambda}
{b_i - b_{i-1}} I)^{-1} \right)\\
&=& C + \sum_{i=1}^{K+1} \psi(b_{i-1},b_i),
\EEAS
where $C=-(Tn/2)(\log(2 \pi) + 1)$
is a constant that does not depend on $b$, and
\[
\psi(b_{i-1},b_i) = -\frac{1}{2} \left(
(b_i-b_{i-1}) \log\det (S^{(i)}+\frac{\lambda}{b_i - b_{i-1}} I)
- \lambda \Tr (S^{(i)}+\frac{\lambda}
{b_i - b_{i-1}} I)^{-1}\right).
\]
(Note that $S^{(i)}$ depends on $b_{i-1}$ and $b_i$.)
Without regularization, \ie, with $\lambda =0$, we have
\[
\psi(b_{i-1}, b_i) =
-\frac{1}{2} (b_i-b_{i-1}) \log\det S^{(i)}.
\]

More generally, we have reduced the regularized maximum likelihood estimation
problem, for fixed values of $K$ and $\lambda$, to the purely combinatorial
problem
\BEQ \label{e-prob}
\begin{array}{ll} 
\mbox{maximize} & -\frac{1}{2}\sum_{i=1}^{K+1} \left(
(b_i-b_{i-1}) 
\log\det (S^{(i)}+\frac{\lambda}{b_i - b_{i-1}} I)
- \lambda \Tr (S^{(i)}+\frac{\lambda}
{b_i - b_{i-1}} I)^{-1} \right),
\end{array}
\EEQ
where the variable to be chosen is the collection of breakpoints 
$b=(b_1, \ldots, b_K)$.  These can take $T-1 \choose K$ possible values.
Note that the breakpoints $b_i$ appear in the objective 
of~(\ref{e-prob}) both 
explicitly and implicitly, through the empirical covariance
matrices $S^{(i)}$, which depend on the breakpoints.

\paragraph{Efficiently computing the objective.}
For future reference, we mention how the objective in \eqref{e-prob} 
can be computed, given $b$.  
We first compute the empirical covariance matrices $S^{(i)}$, which costs
order $Tn^2$ flops.  This step can be carried out in parallel, on
up to $K+1$ processors.  The storage required to store these 
matrices is order $Kn^2$ doubles.
(For comparison, the storage required for the original problem 
data is $Tn$.  Since we typically have $Kn\leq T$, \ie, the average
segment length is at least $n$, the storage of $S^{(i)}$ is no more 
than the storage of the original data.)

For each segment $i=1, \ldots, K+1$, we carry out the following
steps (again, possibly in parallel) to evaluate $\psi(b_{i-1},b_i)$.
We first carry out the Cholesky factorization 
\[
LL^T = S^{(i)}+ \frac{\lambda}{b_i - b_{i-1}} I, 
\]
where $L$ is lower triangular with positive diagonal entries,
which costs order $n^3$ flops.
The log-determinant term can be computed in order $n$ flops,
as $2\sum_{i=1}^{n}\log (L_{ii})$,
and the trace term in order $n^3$ flops, as $\|L^{-1}\|_F^2$.
The overall complexity of evaluating the objective is order
$Tn^2 + Kn^3$ flops, and this can be easily parallelized into $K+1$ independent
tasks.
While we make no assumptions about $T$, $n$, and $K$ (other than $K<T$),
the two terms are equal in order when $T=Kn$, which means that the average segment
length is on the order of $n$, the vector dimension.  This is the threshold at
which the empirical covariance matrices (can) become nonsingular,
though in most applications, useful values of $K$ are much smaller, which means the 
first term dominates (in order).  With the assumption that the average 
segment length is at least $n$, the overall complexity of evaluating
the objective is $Tn^2$.

As an example, we might expect a serial implementation for a data set with
$T=1000$ and $n=100$ to require on the order of 0.01 seconds to 
evaluate the objective, using the very conservative estimate of 1Gflop/sec for
computer speed.

\paragraph{Globally optimal solution.}
The problem \eqref{e-prob} 
can be solved globally by dynamic programming 
\cite{bellman1961approximation, fragkou2004dynamic, kehagias2006dynamic}.
We take as states the set of pairs $(b_{i-1},b_i)$, with $b_{i-1}<b_i$,
so the state space has cardinality $T(T-1)/2$.  We consider the selection 
of a sequence of $K$ states, 
with the state transition constraint that $(p,q)$ must be 
succeeded by a state of the form $(q,r)$.   The complexity of this 
dynamic programming method is $n^3 K T^2$.
Our interest, however, is in a method for large $T$, so we instead seek a
heuristic method for solving \eqref{e-prob} approximately, but with linear
complexity in the time series length $T$.

%The problem \eqref{e-prob} 
%can in principle be solved globally by exhaustive evaluation,
%when $T-1 \choose K$ is not too large (which requires that $K$ be very small, 
%say one or two).
%For other cases, we are not aware of an efficient method for its 
%global solution.  
%Instead, we will seek a heuristic
%method for solving \eqref{e-prob} approximately.
%In general, and also in particular in this application,
%there is little practical benefit (in terms of model
%quality) in globally maximizing a regularized log-likelihood function;
%a good heuristic typically yields models that are just as good when
%applied to real data.

\paragraph{Our method.}
In \S\ref{s-alg}, we describe a heuristic
method for approximately solving problem \eqref{e-prob}.  The method
is not guaranteed to find the globally optimal choice of 
breakpoints, but it does find breakpoints with high (if not always highest)
objective value, and the ones it finds are 1-OPT, meaning that no change of 
any one breakpoint can increase the objective.
The storage requirements of the method are on the order
of the storage required to evaluate the objective, and the computational
cost is typically smaller than a few hundred evaluations of the 
objective function.

\section{Greedy Gaussian segmentation}\label{s-alg}

In this section we describe a greedy algorithm for fitting 
an SGM to data, which we call \emph{greedy Gaussian segmentation} (GGS).
GGS computes an approximate solution of \eqref{e-prob} in a scalable way,
in each iteration adding one breakpoint and then adjusting all the breakpoints
to (approximately) maximize the objective.
In the literature on time series segmentation, this is similar to the standard
``top-down'' approach \cite{keogh2004segmenting}.

\subsection{Split subroutine}

The main building block 
of our algorithm is the Split subroutine. The function $\Split (b_{i-1},b_{i})$ 
takes segment $i$ and finds the $t$ that maximizes 
$\psi (b_{i-1}, t) + \psi (t, b_i)$ over all values of $t$ between $b_{i-1}$ 
and $b_i$. (We assume that $b_i-b_{i-1}>1$; otherwise we cannot split 
the $i$th segment into two segments.)
The time $t=\Split(b_{i-1},b_i)$ is the optimal place to add a breakpoint between
$b_{i-1}$ and $b_i$.
The value of $\psi (b_{i-1}, t) + 
\psi(t, b_i) -\psi (b_{i-1}, b_i)$
is the increase in the objective if we add a new breakpoint at $t$.
This is highest when we choose $t=\Split(b_{i-1},b_i)$.
Due to the regularization term, it is possible for this maximum increase 
to be negative, which means that adding
any breakpoint between $b_{i-1}$ and $b_i$ actually decreases the objective.
The Split subroutine is summarized in Algorithm~\ref{splitAlgo}.

\begin{algorithm}
\caption{Splitting a single interval into two separate segments}
\label{splitAlgo}
\begin{algorithmic}[1]
\Input{$x_{b_{i-1}}, \ldots, x_{b_i}$, along with empirical mean $\mu$
and covariance $\Sigma$.}
\State \textbf{initialize} $\mu_{\mathrm{left}} = 0$, 
$\mu_{\mathrm{right}} = \mu$, $\Sigma_{\mathrm{left}} = \lambda I$, 
$\Sigma_{\mathrm{right}} = \Sigma + \lambda I$.
\For{$t = b_{i-1} + 1, \ldots, b_i - 1$}
\State Update $\mu_{\mathrm{left}}$, $\mu_{\mathrm{right}}$, 
$\Sigma_{\mathrm{left}}$, $\Sigma_{\mathrm{right}}$.
\State Calculate $\psi_t = \psi (b_{i-1}, t)+ \psi (t, b_i)$.
\EndFor
\State\Return The $t$ which maximizes $\psi_t$ and the value of 
$\psi_t - \psi (b_{i-1}, b_i)$ for that $t$.
\end{algorithmic}
\end{algorithm}

In Split, line~3, updating the empirical mean and covariance of the left and right 
segments resulting from adding a breakpoint at $t$, is done in a recursive setting in order $n^2$
flops \cite{welford1962note}.  Line~4, evaluating $\psi_t$, requires order $n^3$ flops, which dominates.
The total cost of running Split is order $(b_i-b_{i-1})n^3$.

\subsection{GGS algorithm}\label{ggsSection}

We can use the Split subroutine to develop a simple greedy 
method for finding good choices of $K$ breakpoints, for $K=1,\ldots, \Kmax$,
by alternating between adding a new break point to the current set of 
breakpoints, and then 
adjusting the positions of all breakpoints until the result is 1-OPT,
\ie, no change of any one breakpoint improves the objective.
This GGS approach is outlined in Algorithm \ref{GreedyAlgo}.

\begin{algorithm}
\caption{Greedy Gaussian segmentation}
\label{GreedyAlgo}
\begin{algorithmic}[1]
\Input{$x_1, \ldots, x_T$, $\Kmax$}.
\State \textbf{initialize} $b_0 = 1$, $b_{1} = T+1$.
\For {$K = 0, \ldots, \Kmax$-1}
\Statex\hskip\algorithmicindent \textbf{\underline{AddNewBreakpoint:}}
\For{$i = 1, \ldots, K+1$}
\State $(t_i, \psi_{\mathrm{increase}})$ = Split$(b_{i-1},b_i)$.
\EndFor 
\If {All $\psi_{\mathrm{increase}}$'s are negative and $K > 0$}
\State \textbf{return} $(b_1, \ldots, b_K)$.
\ElsIf {All $\psi_{\mathrm{increase}}$'s are negative}
\State \textbf{return} $()$.
\EndIf
\State Add a new breakpoint at the $t_i$ with the largest corresponding 
value of $\psi_{\mathrm{increase}}$.
\State Relabel the breakpoints so that 
$1 = b_0 <b_1< \cdots < b_{K+1} < b_{K+2} = T+1$.
\Statex\hskip\algorithmicindent \textbf{\underline{AdjustBreakpoints:}}
\Repeat
\For{$i = 1, \ldots, K$}
\State $(t_i, \ell_{\mathrm{increase}})$ = Split$(b_{i-1}, b_{i+1})$.
\State If $t_i \neq b_i$, set $b_i = t_i$.
\EndFor
\Until Stationary.
\EndFor
\State \Return  $(b_1, \ldots, b_K)$.
\end{algorithmic}
\end{algorithm}

In line~2, we loop over the addition of new breakpoints, adding
exactly one new breakpoint each iteration. %so that 
%after iteration $K$ we have exactly $K$ breakpoints.
Thus, the algorithm finds good sets of breakpoints, for
$K=1, \ldots, \Kmax$, unless it quits early in
line~6.  This occurs when the addition of any new breakpoint will 
decrease the objective.
In AdjustBreakpoints, we loop over the current segmentation 
and adjust each 
breakpoint alone to maximize the objective. In this step the objective can
either increase or stay the same, and we repeat until the current
choice of breakpoints is 1-OPT. 
In AdjustBreakpoints, there is no need to call 
Split$(b_{i-1}, b_{i+1})$ more than once 
if the arguments have not changed.

The outer loop over $K$ must be run 
serially, since in each iteration we start with the breakpoints 
from the previous iteration.
Lines~3 and~4 (in AddNewBreakpoint) can be run in parallel 
over the $K+1$ segments.
We can also parallelize AdjustBreakpoints,
by alternately adjusting the even and odd breakpoints (each 
of which can be parallelized) until stationarity.
GGS requires storage on the order of $Kn^2$ numbers.
As already mentioned, this is typically the same order as, or less
than, the storage required for the original data.

Ignoring opportunities for parallelization,
running iteration $K$ of GGS requires order
$KLn^3T$ flops, where $L$ is the average number of iterations required in
AdjustBreakpoints.
When parallelized, the complexity drops to $Ln^3T$ flops.
While we do not know an upper bound on $L$, 
we have observed empirically that it is modest when $K$ is not too large;
that is, AdjustBreakpoints runs just a few outer loops over the breakpoints.
Summing from $K=1$ to $K=\Kmax$, and assuming $L$ is a constant,
gives a complexity of order
$(\Kmax)^2 n^3 T$ (without parallelization), or 
$\Kmax n^3T$ (with parallelization).
In contrast, the dynamic programming method
\cite{bellman1961approximation, fragkou2004dynamic, kehagias2006dynamic}
requires order $\Kmax n^3 T^2 $ flops.

\section{Validation and parameter selection}\label{s-validation}

Our GGS method has just two hyperparameters: 
$\lambda$, which controls the 
amount of (inverse) covariance regularization, and $\Kmax$, the 
maximum number of breakpoints.
In applications where the reason for segmentation is to
identify interesting times where the statistics of 
the data change, $K$ (and $\lambda$) might be chosen by hand, or
by aesthetic or other considerations, such as whether the segmentation
identifies known or suspected times when something changed.
The hyperparameter values can also be chosen by a more principled method,
such as Bayesian or Akaike information criterion \cite[chapter 7]{hastie2009elements}.
In this section, we describe a simple method of 
selecting the hyperparameters
through out-of-sample or cross validation.
We first describe the basic idea with 10:1 out-of-sample validation.

We remove 10\% of the data at random,
leaving us with $0.9T$ remaining samples.  The 10\% of samples
are our test set, and the remaining samples are the training set,
which we use to fit our model.  We choose some reasonable value
for $\Kmax$, such as $\Kmax= (T/n)/3$ (which corresponds to the average 
segment length $3n$) or a much smaller number when $T/n$ is large.
For multiple values of $\lambda$, typically logarithmically
spaced over a wide range, we run the GGS algorithm.
This gives us one SGM for each value of $\lambda$ and each value of $K$.
For each of these SGMs, we note the log-likelihood on
the training data, and also on the test data.
(It is convenient to divide each of these by the number of data 
points, so they become the average log-likelihood per sample. 
In this way the numbers for the train and test sets can be compared.)
%We can speed up this computation by running GGS 
%in parallel for each value of $\lambda$.
To calculate the log-likelihood on the test set, we simply evaluate
\[
\ell(x_t) =
-\frac{1}{2} (x_t - \mu^{(i)})^T (\Sigma^{(i)})^{-1} (x_t - \mu^{(i)}) - 
\frac{1}{2} \log\det \Sigma^{(i)} - \frac{n}{2}\log(2\pi),
\]
if $t$ falls in the $i$th segment of the model.
The overall test set log-likelihood is then defined, on a test set
$\mathcal{X}$, as 
\[
\frac{1}{| \mathcal{X} |}\sum_{x_t \in \mathcal{X}} \ell(x_t).
\]
Note that when $t$ is the time index of a sample in the test set,
it cannot be a breakpoint of the model, since the model was developed 
using the data in the training set.

We then apply standard principles of validation.  If for a particular 
SGM (found by GGS with a particular value of $\lambda$ and $K$) 
the average log-likelihood on the train and test sets is similar,
we conclude the model is not over-fit, and therefore 
a reasonable candidate.
Among candidate models, we then choose one that has a high value of 
average log-likelihood.  If many models have reasonably high 
average log-likelihood, we choose one with a small value of $K$ and 
a large value of $\lambda$.  (In the former case to get the simplest 
model that explains the data, and in the latter case to get the least
sensitive model that explains the data.)

Standard cross-validation is an extension of out-of-sample validation 
that can give us even more confidence in a proposed SGM.
In cross-validation we divide the original data into 10 equal size `folds'
of randomly chosen samples, and carry out out-of-sample validation
10 times, with each fold as the test set.
If the results are reasonably consistent across the folds, both in
terms of train and test average log-likelihood and the breakpoints
themselves, we can have confidence that the SGM fits the data.

\section{Variations and extensions}\label{s-variations}

The basic model and method can be extended in many ways,
several of which we describe here.

\paragraph{Warm-start.}
GGS builds SGMs by increasing $K$, starting from $K=0$.
It can also be used in warm-start mode, meaning we start the algorithm
from a given choice of initial breakpoints.
As an extreme version, we can start with a random set of $K$ breakpoints,
and then run AdjustBreakpoints until we have a 1-OPT solution. 
The main benefit of a warm start is that it allows for a significant computational
speedup. Whereas a (parallelized) GGS algorithm has a runtime of $O(\Kmax Tn^3)$,
this warm-start method takes only $O(Tn^3)$, since it can skip the first
$\Kmax - 1$ steps of Algorithm 2. 
However, as we will show in~\S 6.2, this speedup comes with a tradeoff, as the 
solution accuracy tends to drop when running GGS in warm-start mode as compared
to the original algorithm.

\paragraph{Backtracking.}
In GGS we add one breakpoint per iteration.  
While we adjust the previous breakpoints found, we never remove a breakpoint.
One variation is to occasionally remove a breakpoint.  This can be done 
using a subroutine called Combine.
This function evaluates, for each breakpoint, the decrease in objective
value if that breakpoint is removed.
In a backtracking step, we remove the breakpoint that decreases the objective the
least; we then can adjust the remaining breakpoints, and continue with
the GGS algorithm, adding a new breakpoint.  (If we end up adding the breakpoint
we removed back in, nothing has been achieved.) 
We also note that backtracking allows for GGS to be solved by a bottom-up method 
\cite{keogh2004segmenting, borenstein2008combined}. We do so by starting with 
$T-1$ breakpoints and continually backtracking until only $K$ breakpoints remain.

%\begin{figure}
%\caption{Backtracking algorithm can help the GGS algorithm avoid suboptimal solutions.}
%\centering
%\includegraphics[width=0.6\textwidth]{figures/backtracking.pdf}
%\label{fig:backtracking}
%\end{figure}

\paragraph{Streaming.}
We can deploy GGS when the data is streaming.
We maintain a memory of the last $M$ samples and run GGS on this data set.
We could do this from scratch as each new data point or group of data points arrives,
complete with selection of the hyperparameters and validation.  
Another option is to fix $\lambda$ and $K$, and then to run GGS in warm-start mode, which means
that we keep the previous breakpoints (shifted appropriately), and then run 
AdjustBreakpoints from this starting point (as well as AddBreakpoint if a breakpoint
has fallen off our memory).

In streaming mode, the GGS algorithm provides an estimate of the statistics of 
future time samples, namely, the mean and covariance in the SGM in the most recent segment.
%XXX this is important and we should try this.

%look for a fast way to incorporate the $T+1$st reading.
%To do so, we first run Split$(b_k, T+1)$ and add the $K+1$st
%breakpoint there. Then, we run the Combine function, seeing how much
%%the likelihood would decrease by removing each of the $K+1$ breakpoints
%(we can cache the results for the first $K-1$ $b_i$'s). If Combine
%returns the $K+1$st breakpoint, then we remove it and continue. However,
%if it returns another breakpoint, we determine that the timestamp returned
%by Split$(b_k, T+1)$ is a legitimate breakpoint, and we keep it.

\paragraph{Multiple samples at the same time.}
Our approach can easily incorporate the case where we have more
than one data vector for any given time $t$.  
We simply change the sums over each segment for the empirical mean and 
covariance to include any data samples in the given time range.
%have multiple vector readings $x$ at a single time $t$. And so, we may
%have $W$ samples spread out over $T$ different timestamps. In this setup,
%we can run the GGS algorithm as usual over the $T$ timestamps, and the 
%only difference is in the Split function, on Step 3 of Algorithm \ref{splitAlgo}.
%When we update the means and covariances of the left and right split, we
%simply incorporate \emph{all} samples at each time $t$. The runtime of
%the algorithm remains nearly identical, going from $Tn^3 + Tn^2$ for the 
%original algorithm to $Tn^3 + Wn^2$ for this new one.

\paragraph{Cyclic data.}
In cyclic data, the times $t$ are interpreted modulo $T$, so $x_T$ and $x_1$
are adjacent.
A good example is a vector time series that represents daily measurements over
multiple years; we simply map all measurements to $t=1,\ldots, 365$ (ignoring
leap years), and modify the model and method to be cyclic.
The only subtlety here arises in choosing the first breakpoint,
since one breakpoint does not split a cyclic set of times into 
two segments. Evidently we need two breakpoints to split a cyclic 
set of times into two segments.
We modify GGS by arbitrarily choosing a first breakpoint, 
and then running as usual, including the `wrap-around' segment as a segment.
Thus, the first step chooses the second breakpoint, which splits the cyclic
data into two segments.  The AdjustBreakpoints method now adjusts both the 
chosen breakpoint and the arbitrarily chosen original one.
%XXX we should implement this and try it.

\paragraph{Regularization across time.}
In our current model, the estimates on either side of a breakpoint are independent
of each other.
We can, however, carry out a post-processing step to shrink
models on either side of each breakpoint towards each other.
We can do this by fixing the breakpoints and then adjusting the 
continuous model parameters to minimize our original objective minus
a regularization term that penalizes deviations 
of $(\Sigma^{(i)},\mu^{(i)})$ from $(\Sigma^{(i-1)},\mu^{(i-1)})$.
%\[
%\sum_{i=1}^{K+1} \alpha 
%\left\|(\Sigma^{(i)} + \lambda I)^{-1} - (\Sigma^{-1(i-1)} \|_2 + 
%\beta \|\Sigma^{-1(i)}\mu^{(i)} - \Sigma^{-1(i-1)}\mu^{(i-1)} \|_2.
%\]
%(This choice of regularization yields a convex problem in
%the natural variables $\Sigma^{-1}$ and $\Sigma^{-1}\mu$.)

%This is equivalent to $\ell_2$ temporal regularization, rather than
%the $\ell_0$ regularization that our original model assumes.
%Note that this couples the $\mu$'s and $\Sigma$'s together, 
%so we need to jointly solve for them all, which we can do by solving
%for the covariance-regularized likelihood plus temporal regularization
%problem for $\mu$ and $\Sigma$, which now more feasible, 
%since $K \ll T$.

\paragraph{Non-Gaussian data.}
Our segmented Gaussian model and associated regularized maximum likelihood 
problem \eqref{e-prob} can be generalized to other statistical models.
The problem is tractable, at least in theory, when the associated regularized
maximum likelihood problem is convex.  In this case we can compute the optimal
parameters over a segment by solving a convex optimization problem, whereas in the SGM
we have an analytical solution in terms of the empirical mean and covariance.
Thus we can segment Poisson or Bernoulli data, or even
heterogeneous exponential family distributions 
\cite{lee2015learning, tansey2015vector}.
%when some elements of the vector 
%come from a Poisson distribution while others are Bernoulli.

\section{Experiments}\label{s-examples}

In this section, we describe our implementation of GGS, and the results of some 
numerical experiments to illustrate the model and the method.

\subsection{Implementation}
We have implemented GGS as a Python package \texttt{GGS} available at 
\begin{center}
\url{https://github.com/cvxgrp/GGS}.
\end{center}
\texttt{GGS} is capable of carrying out full ten-fold cross-validation 
to help users choose values of the hyperparameters.
\texttt{GGS} uses \texttt{NumPy}  
for the numerical computations, and
the \texttt{multiprocessing} package 
to carry out the algorithm in
parallel for different cross-validation folds for a single $\lambda$.
(The current implementation does not support parallelism over the segments
of a single fold, and the advantages of parallelism will only be seen
when \texttt{GGS} is run on a computer with multiple cores.)

\subsection{Financial indices}

In financial markets, regime changes have been shown to have
important implications for asset class and portfolio performance
\cite{ang2012regime, sheikh2012regime, nystrup2015regime, nystrup2017dynamic}.
We start with a small example with $n=3$, 
where we can visualize and plot all entries of the segment 
parameters $\mu^{(i)}$ and $\Sigma^{(i)}$.

\begin{figure}
\caption{Cumulative returns over the 19-year period for a stock, oil, and bond index.}
\centering
\vspace{-3mm}
\includegraphics[width=0.7\textwidth]{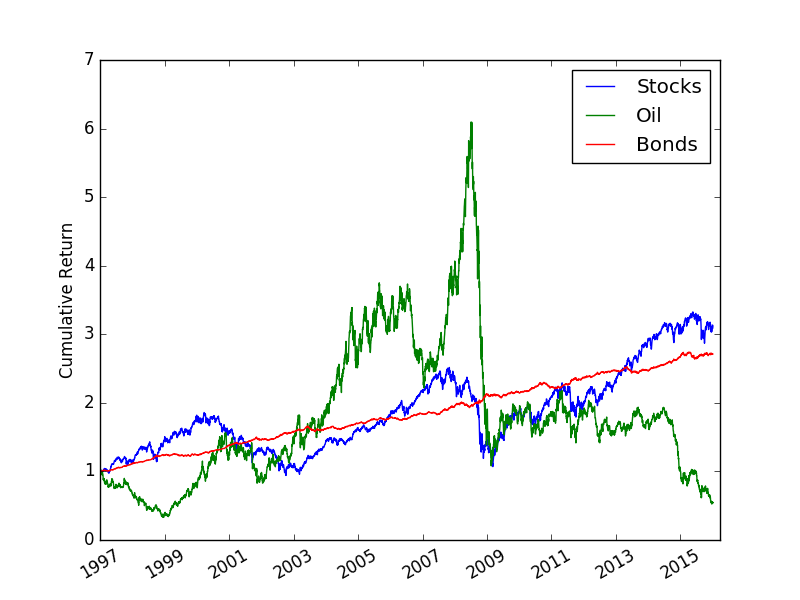}
\label{fig:cumulative}
\end{figure}

\paragraph{Data set description.}
Our data set consists of 19 years of daily returns, from January 1997 to 
December 2015, for $n=3$ indices for stocks, oil, and 
government bonds: MSCI World, S\&P GSCI
Crude Oil, and J.P.\ Morgan Global Government
Bonds. 
We use log-return data, \ie, the logarithm of the end-of-day price increase
from the previous day.
The time series length is $T=4943$.
Cumulative returns for the three indices are shown in Figure \ref{fig:cumulative}.
We can clearly see multiple `regimes' in the return series, although the 
individual behaviors of the three indices are quite different.

\begin{figure}
\caption{Objective $\phi(b)$ vs. number of breakpoints $K$ for $\lambda =10^{-4}$.}
\centering
\includegraphics[width=0.7\textwidth]{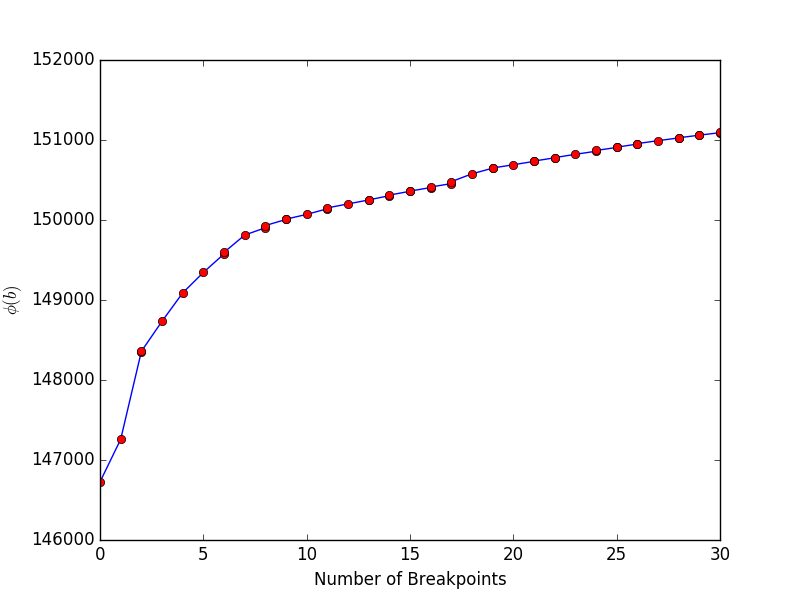}
\label{fig:regpath}
\end{figure}

\paragraph{GGS algorithm.}
We run GGS on the data with $\Kmax=30$ and $\lambda=10^{-4}$.
Figure \ref{fig:regpath} shows the objective function versus $K$,
\ie, the objective in each iteration of GGS.
We see a sharp increase in the objective up to around $K = 8$ or $K=10$ --- 
our first hint that a choice in this range would be reasonable.
For this example $n$ is very small, so the computation time is dominated 
by Python overhead.  Still, our single-threaded GGS solver took less than 
30 seconds to compute these 30 models on a standard laptop with a
1.7 GHz Intel i7 processor.
The average number of passes through
the data for the breakpoint adjustments was under two.

%\paragraph{Out-of-sample validation.}
%We reserve 10\% of the data points at random as a test set, and 
%fit GGMs using GGS with the remaining training data with regularization
%parameter values $\lambda = 0.01, ~0.1,~1.0,~10$, and $\Kmax = 30$.
%The results are shown in Figure \ref{fig:testSets}.
%We can see that for small and even medium regularization, the average log-likelihood 
%in the test and training sets deviate for $K$ larger than around 8 or so,
%meaning the model is over-fit.
%Among these, the highest test likelihood is obtained for $\lambda =1$ (although
%$\lambda = 0.1$ also does well), and a value of $K$ up to around 10
%would be reasonable.

\begin{figure}
\caption{Average train and test set log-likelihood during 10-fold cross validation
for various $\lambda$'s and across all values of $K \leq 30$.}
\begin{subfigure}{0.5\textwidth}
  \centering
  \includegraphics[width=.99\linewidth]{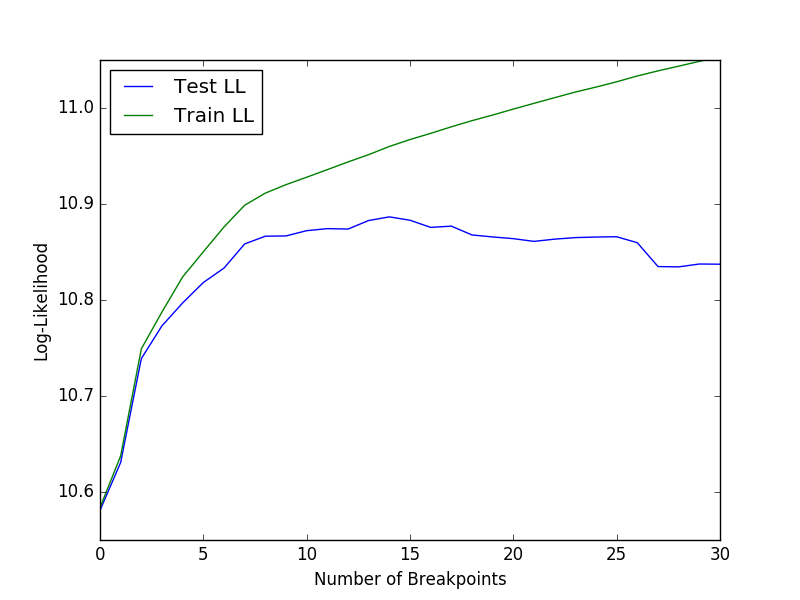}
  \caption{$\lambda = 10^{-6}$}
  \label{fig:sfig1}
\end{subfigure}%
\begin{subfigure}{.5\textwidth}
  \centering
  \includegraphics[width=.99\linewidth]{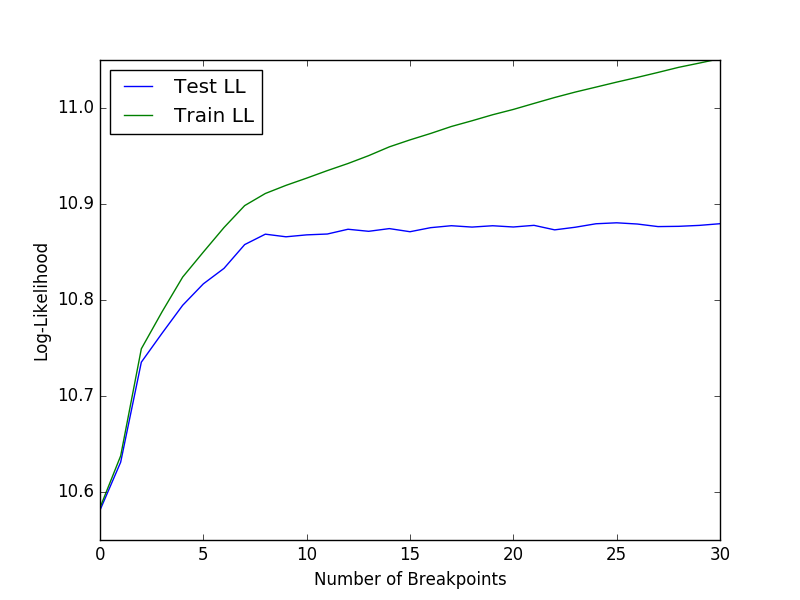}
  \caption{$\lambda = 10^{-5}$}
  \label{fig:sfig2}
\end{subfigure}
\begin{subfigure}{.5\textwidth}
  \centering
  \includegraphics[width=.99\linewidth]{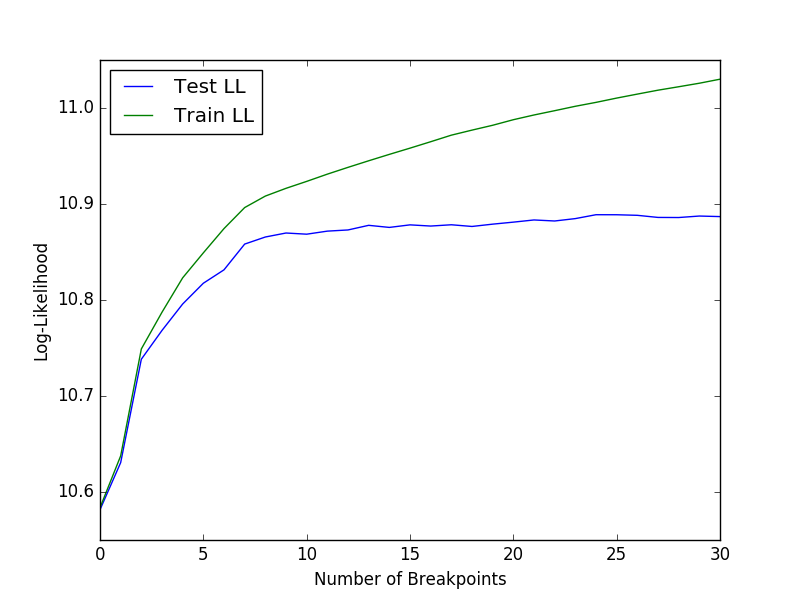}
  \caption{$\lambda = 10^{-4}$}
  \label{fig:sfig3}
\end{subfigure}%
\begin{subfigure}{.5\textwidth}
  \centering
  \includegraphics[width=.99\linewidth]{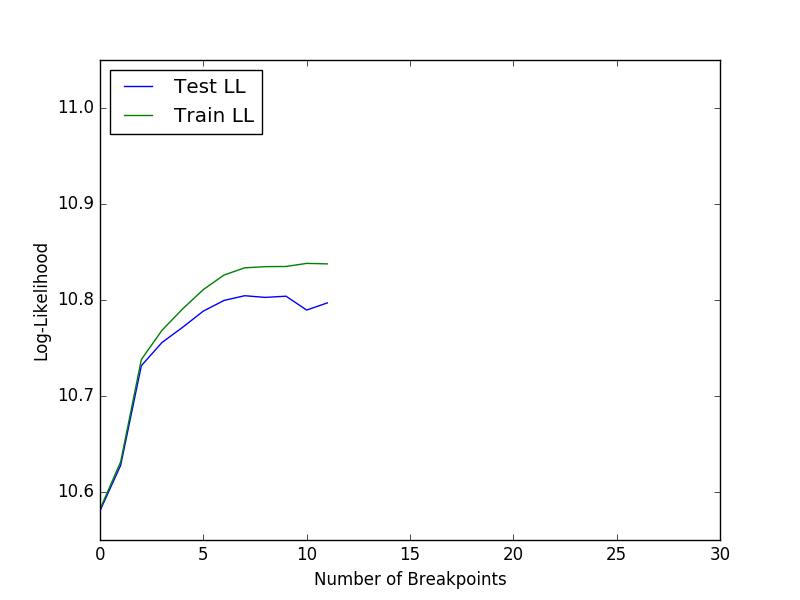}
  \caption{$\lambda = 10^{-3}$}
  \label{fig:sfig4}
\end{subfigure}
\label{fig:validation}
\end{figure}

\paragraph{Cross-validation.}
We next use 10-fold cross-validation to determine reasonable values 
for $K$ and $\lambda$. We plot the average log-likelihood over the 10 folds
versus $K$ in Figure \ref{fig:validation} for various values of 
$\lambda$. 
When $\lambda$ is large, the curves stop before $K = \Kmax$, since GGS 
terminates early.
These plots clearly show that increasing $K$ above 10 does not 
increase the average log-likelihood in the test set;
and moreover past this point the
log-likelihood on the test and training sets begin to diverge,
meaning the model is overfit.
Though Figure \ref{fig:validation} only goes up to $K=30$, we find
that for values of $K$ above around 60,
the log-likelihood begins to drop significantly.
Furthermore, we see that values of 
$\lambda$ up to $\lambda = 10^{-4}$ yield roughly the same high log-likelihood.
This suggests that choices of $K = 10$ and 
$\lambda = 10^{-4}$ are reasonable, aligning with our general
preference for models which are simple (small $K$) and not too sensitive 
to noise (large $\lambda$).
Cross-validation also reveals that the choice of breakpoint locations is 
very stable for these values of $K$ and $\lambda$, across the 10 folds.

%\begin{figure}
%\caption{Cross validation results, with vertical axis showing 
%average log-likeligood on train and test sets.}
%\centering
%\includegraphics[width=0.7\textwidth]{figures/testSets.png}
%\label{fig:testSets}
%\end{figure}

\begin{figure}
\caption{Segmented Gaussian model obtained with
$\lambda = 10^{-4}$, $K = 10$.}
\centering
\includegraphics[width=0.95\textwidth]{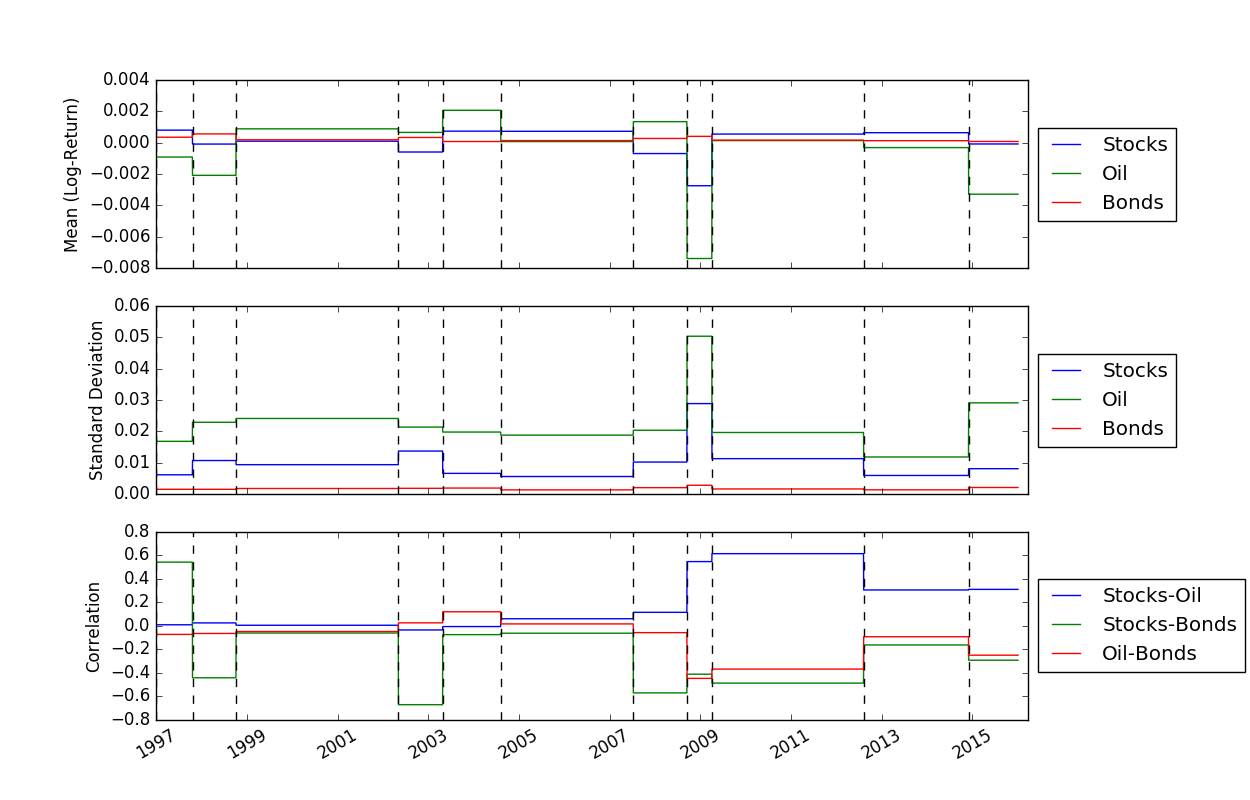}
\label{fig:timeSeries}
\end{figure}

\begin{figure}
\caption{Cumulative returns with vertical bars at the model breakpoints.}
\centering
\includegraphics[width=0.75\textwidth]{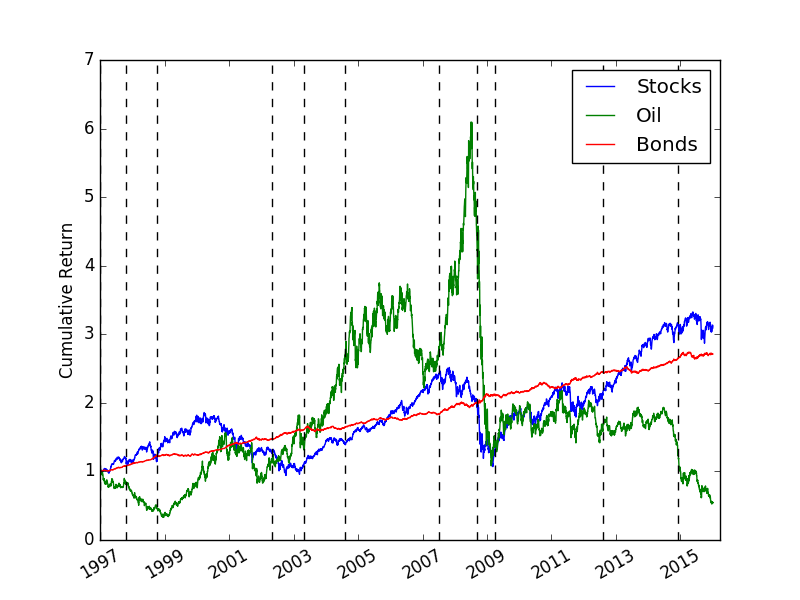}
\label{fig:cumulativeVBars}
\end{figure}

\paragraph{Results.}
Figure \ref{fig:timeSeries} shows the model obtained by GGS with 
$\lambda =10^{-4}$ and $K=10$.
We plot the covariance matrix by showing the square root of
the diagonal entries (\ie, the volatilities), and the three correlations,
versus $t$.
During the financial crisis in 2008, the mean return of stocks and oil was
very negative and the volatility was high. 
The stock market and the oil price were almost uncorrelated before 2008,
but have been positively correlated since then. It is interesting to see
how the correlation between stocks and bonds has varied over time: It was
strongly positive in 1997 and very negative in 1998, in 2002, and in the five years
from mid-2007 to mid-2012. 
The sudden shift in this correlation between 1997 and 1998 is why GGS yields two
relatively short segments in the [1997, 1999] window, rather than breaking up a longer
segment (such as [1999, 2002], where the correlation structure is more homogenous).
The extent of these variations would be difficult to capture using
a sliding window; the window would have to be very short, which would
lead to noisy estimates. The segmentation approach yields a more interpretable
partitioning with no dependence on a (prespecified, fixed) window length.
Approaches to risk modeling \cite{alexander2000primer} and portfolio optimization
\cite{partovi2004principal, meucci2009managing} based on principal
component analysis are questionable, when volatilities and correlations are
changing as significantly as is the case in Figure \ref{fig:timeSeries}
(see also \cite{fenn2011temporal}).
We plot the cumulative index returns along with the chosen breakpoints in
Figure \ref{fig:cumulativeVBars}.
We can clearly see natural segments and boundaries, for example the 
Russian default in 1998 and the 2008 financial crisis.

\begin{figure}
\caption{Empirical complementary CDF of $\phi(b)$ for 10000 randomly initialized
results for $K = 10$, $\lambda = 10^{-4}$.
Vertical bars represent GGS solutions for different values of $K$.} %$K = 8$ (red), $K = 9$ 
%(green), $K = 10$ (black), and $K = 11$ (purple).}
\centering
\includegraphics[width=0.75\textwidth]{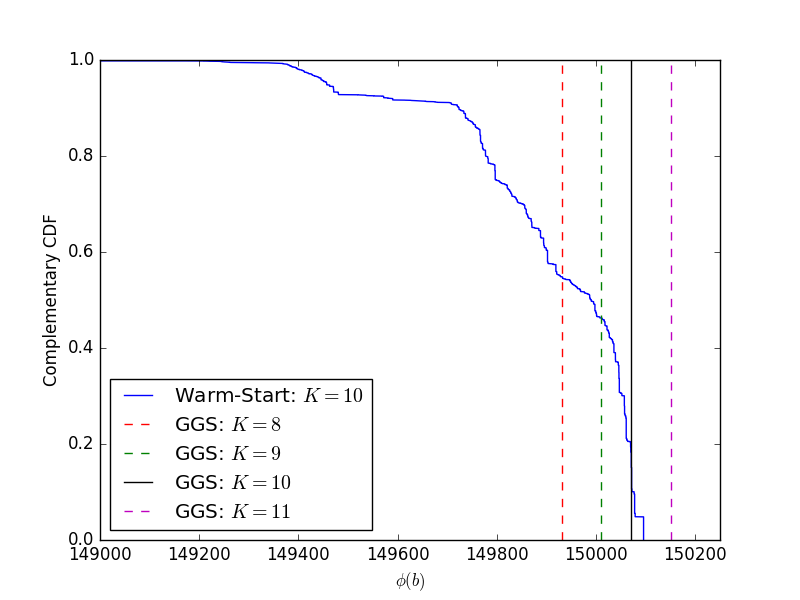}

\label{fig:ccdf}
\end{figure}

\paragraph{Comparison with random warm-start.}
We fix the hyperparameters $K = 10$ and $\lambda = 10^{-4}$,
and attempt to find a better SGM using warm-start with 
random breakpoints.
This step is not needed; we carry this out to demonstrate that 
while GGS does not find the model that globally maximizes the objective,
it is effective.
We run 10000 warm-start random initial breakpoint computations,
running AdjustBreakpoints until the model is 1-OPT
and computing the objective found in each case.
(In this case the number of passes over the data set far exceeds two, 
the typical number in GGS.)
The complementary CDF of the objective for these 10000
computations is shown in Figure \ref{fig:ccdf},
as well as the objective values found by GGS
for $K = 8$ through $K = 11$. 
We see that the random initializations can sometimes 
lead to very poor results; over 50\% of the simulations,
even though they are locally optimal,
have smaller objectives than the $K=9$ step of GGS. 
On the other hand, the random initializations do find some 
SGMs with objective slightly exceeding the one found by GGS,
demonstrating that GGS did not find the globally optimal set of breakpoints.
These SGMs have similar breakpoints, and similar cross-validated log-likelihood,
as the one found by GGS. As a practical matter, these SGMs are no better
than the one found by GGS.
There are two advantages of GGS over the random search: First, it is much faster;
and second, it finds models for a range of values of $K$, which is useful before
we select the value of $K$ to use.

\clearpage
\subsection{Large-scale financial example}

\paragraph{Data set description.}
We next look at a larger example to emphasize the scalability 
of GGS. We look at all companies currently in the S\&P 
500 index that have been publicly listed for the entire 19-year interval from 
before (from 1997 to 2015), which leaves 309 companies. 
Note that there are slightly fewer trading days for the S\&P 500
each year than the global indices, since the S\&P 500 does not
trade during US holidays, while the global indices still move. 
The 19-year data set yields a $309 \times 4782$ data matrix.
We take daily log-returns for these stocks, and run the GGS 
algorithm to detect relevant breakpoints.

\paragraph{GGS scalability.}
We run GGS on this much larger data set up to $\Kmax=10$. 
%in Figure 
%\ref{fig:bigGGS}. (TODO -remove this plot?)
Our serial implementation of the GGS algorithm, on the same 1.7 GHz 
Intel i7 processor, took 36 minutes, where 
AdjustBreakpoint took an average of 3.5 passes through the data
at each $K$. Note that this aligns very closely with our predicted
runtime from \S\ref{ggsSection}, which was estimated as $(\Kmax)^2 L T n$.

%2185.67938304 seconds

%\begin{figure}
%\caption{$\phi(b)$ vs. the number of breakpoints for the entire 19-year
%interval for the S\&P 500 for $\lambda = 10^{-2}$.}
%\centering
%\includegraphics[width=0.6\textwidth]{figures/LLbigExample.png}
%\label{fig:bigGGS}
%\end{figure}

\begin{figure}
\caption{Average train and test set log-likelihood during 10-fold cross validation
for various $\lambda$'s for the 309-stock example. Note that not all $\lambda$'s
go all the way up to $K = 10$ because our algorithm stops
when it determines that it will no longer benefit from adding an additional split.}
\begin{subfigure}{0.5\textwidth}
  \centering
  \includegraphics[width=.99\linewidth]{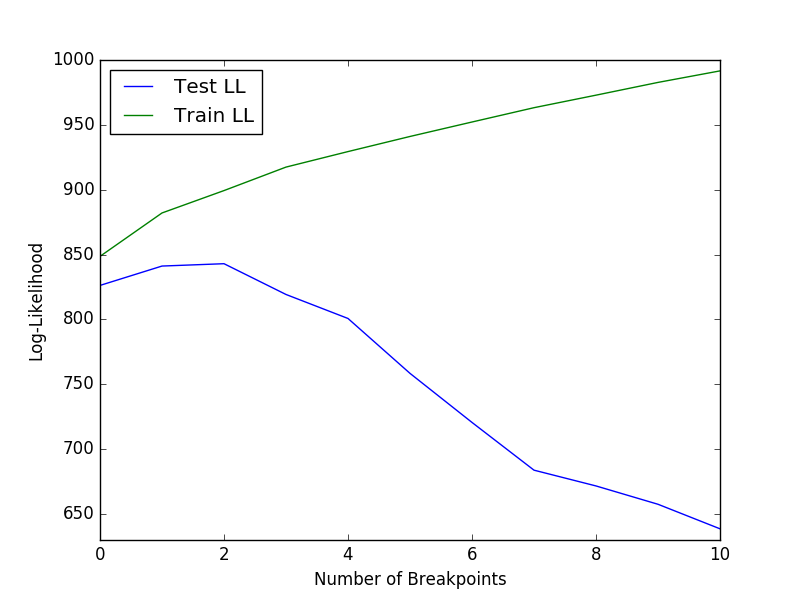}
  \caption{$\lambda = 5\times10^{-3}$}
  \label{fig:sfig1}
\end{subfigure}%
\begin{subfigure}{.5\textwidth}
  \centering
  \includegraphics[width=.99\linewidth]{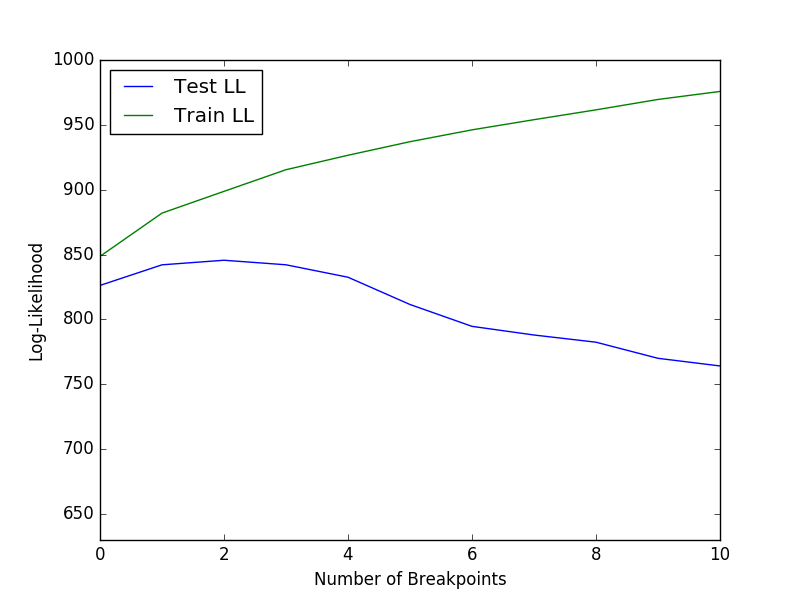}
  \caption{$\lambda = 10^{-2}$}
  \label{fig:sfig2}
\end{subfigure}
\begin{subfigure}{.5\textwidth}
  \centering
  \includegraphics[width=.99\linewidth]{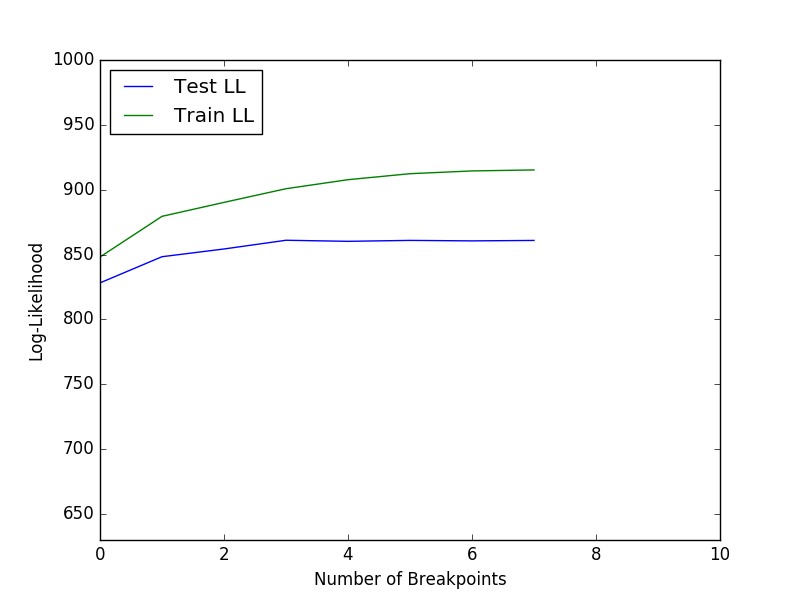}
  \caption{$\lambda = 5\times10^{-2}$}
  \label{fig:sfig3}
\end{subfigure}%
\begin{subfigure}{.5\textwidth}
  \centering
  \includegraphics[width=.99\linewidth]{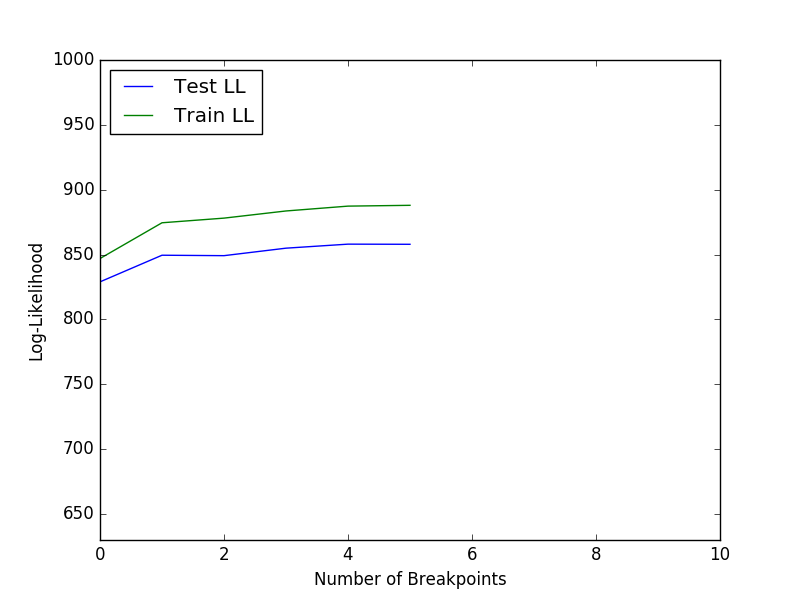}
  \caption{$\lambda = 10^{-1}$}
  \label{fig:sfig4}
\end{subfigure}
\label{fig:bigVal}
\end{figure}

\paragraph{Cross-validation.}
We run 10-fold cross-validation to find good values of the hyperparameters 
$K$ and $\lambda$. The average log-likelihood of the test and train sets
are displayed in Figure \ref{fig:bigVal}. From the results, we can see that
the log-likelihood is maximized at a much smaller value of $K$, indicating
fewer breakpoints. This is in part because, with $n = 309$, we need more
samples in each segment to get an accurate estimate of the $309\times309$ 
covariance matrix, as opposed to the $3\times3$ covariance in the smaller example.
%This is likely because each time segment needs to estimate 
%a 309$\times$309 covariance matrix, which requires many more samples than
%a simple 3$\times$3 covariance, so the maximum likelihood is more conservative
%with its estimates of when a breakpoint as occurred.
Our cross-validation results suggest choosing $K=3$ and $\lambda=5\times10^{-2}$,
and as in the small example, the results are very stable near these values.

\begin{figure}
\caption{Mean, standard deviation, and cumulative return for a uniform
portfolio with $K = 3$, $\lambda = 5\times10^{-2}$.}
\centering
\includegraphics[width=0.99\textwidth]{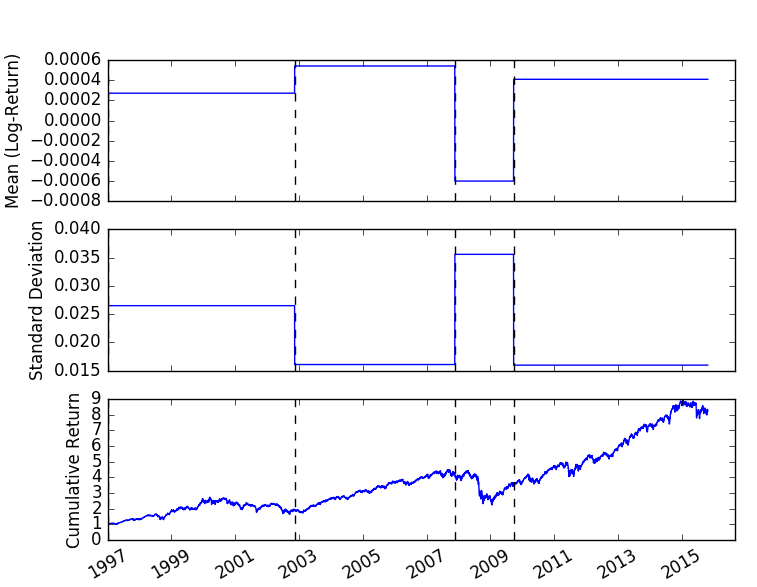}
\label{fig:bigMean}
\end{figure}

\paragraph{Results.}
We plot the mean, standard deviation, and cumulative return of a
uniform, buy-and-hold portfolio (\ie, investing \$1 into each of the 309 stocks 
in 1997). The results are shown in 
Figure \ref{fig:bigMean}. Note that there is a selection bias
in the data set, since these companies all remained in the S\&P 500 
in 2016, and the total return is 8x over the 19-year interval. Like 
before, the 2008 financial crisis stands out. It is the only segment with
a negative mean value. The partitioning seems intuitively right.
The first, highly volatile, segment includes both the build-up and burst
of the dot--com bubble. The second segment is the bull market that led
to the financial crisis in 2008. The third segment is the financial crisis and the
fourth segment is the market rally that followed the crisis. These break points were
also found in the multiasset example in
Figure \ref{fig:timeSeries} and \ref{fig:cumulativeVBars}. 

\clearpage
\subsection{Wikipedia text data}

We examine an example from the field of natural language processing (NLP) to
illustrate how GGS can be applied to a very different type of data set, beyond
traditional time series examples.

\paragraph{Data set description.} We look at text data from English-language
Wikipedia. We obtain our data by concatenating the introductions, \ie, the
text that precedes the Table of Contents section on the Wikipedia webpage,
of five separate articles, with titles George Clooney, Botany, Julius Caesar, 
Jazz, and Denmark. 
Here, the ``time series'' consists of the sequence of 
words from these five articles in order.
After basic preprocessing (removing words that 
appear at least five times and in multiple articles), 
our data set consists of 1282 words, with each article
contributing between 224 and 286 words. We then convert each word into a
300-dimensional vector using a pretrained Word2Vec embedding of
three million unique words (or short phrases), trained on the Google News data set
of approximately 100 billion words, available at
\begin{center}
\url{https://code.google.com/archive/p/word2vec/}.
\end{center}
This leaves us with
a $300 \times 1282$ data matrix.
Our hope is that GGS can detect the breakpoints between the five concatenated 
articles, based solely on the change in mean and covariance of the vectors 
associated with the words in our vector series.

\begin{figure}
\caption{Actual and GGS predicted breakpoints for the concatenation of the 
five Wikipedia articles, along with the predicted most similar word to the mean of each GGS segment.}
\centering
\includegraphics[width=0.99\textwidth]{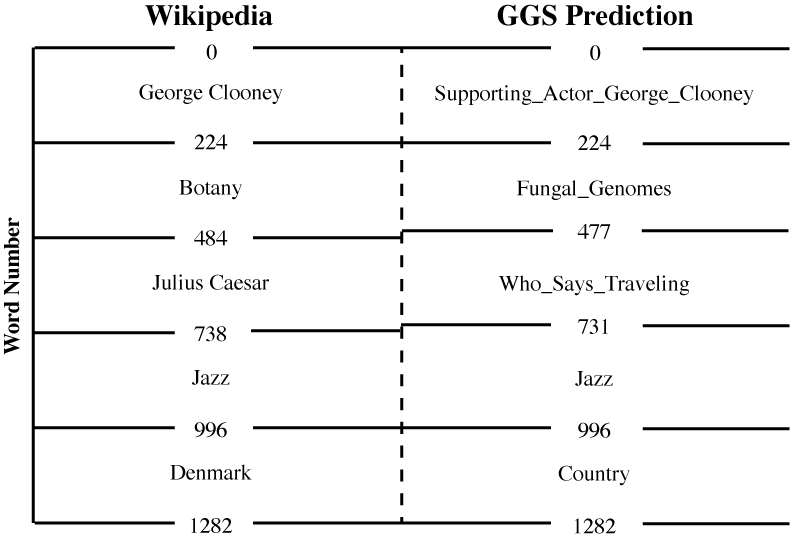}
\label{fig:NLP}
\end{figure}

\paragraph{GGS results.} We run GGS to split the data into five segments --- \ie, $K = 4$ --- and
use cross-validation to select $\lambda = 10^{-3}$.
(We note, however, that this example is quite robust to the selection of $\lambda$, and any 
value from $10^{-6}$ to $10^{-3}$ yields the exact same breakpoint locations.)
We plot the results in Figure \ref{fig:NLP}, which show
GGS achieving a near-perfect split of the five articles. Figure
\ref{fig:NLP} also shows a representative word (or short phrase) in the Google
News data  set that is among the top five ``most similar'' words, out of the entire
three million word corpus, to the average (mean) of each GGS segment, as measured
by cosine similarity. We see that GGS correctly identifies both the
breakpoint locations and the general topic of each segment. 
%This demonstrates
%how GGS can be used for topic modeling of text data, for example discovering
%when a talk radio show has changed discussion topics. Or, it could be used to
%help an editor split a book into ``sections'' to help potential readers
%navigate the book.

\clearpage
\subsection{Comparison with left-to-right HMM on synthetic data}

We next analyze a synthetic example where observations are generated 
from a given sequence of segments. This provides a known ground 
truth, allowing us to compare GGS with a common baseline, a 
left-to-right hidden Markov model (HMM)
\cite{bakis1976continuous, cappe2005inference}. 
Left-to-right HMMs, like GGS, split the data into
non-repeatable segments, where each segment is defined by a Gaussian
distribution. The HMMs in this experiment are implemented using the \texttt{rarhsmm}
library \cite{rarhsmm}, which includes the same shrinkage estimator for the
covariance matrices used in GGS (see \S\ref{s-regularized}).
The shrinkage estimator results in more reliable estimates not only of the covariance
matrices but also of the transition matrix and the hidden states \cite{fiecas2017shrinkage}.

\paragraph{Data set description.} We start by generating 10 random covariance 
matrices. We do so by setting $\Sigma_i = A^{(i)} A^{(i)T}, i = 1, \ldots, 10$, where 
$A^{(i)} \in \reals^{25 \times 25}$ is a random matrix where each element $A^{(i)}_{j,k}$ was 
generated independently from the standard normal distribution.
Our synthetic data set then has 10 ground truth segments (or $K = 9$ breakpoints), where 
segment $i$ has zero mean and covariance $\Sigma_i$. Each segment is of length 100 (so the
total time series has $T=1000$ observations). Each of the 100 readings per segment
is sampled independently from the given distribution. Thus, our final data set consists of a
25 $\times$ 1000 data matrix, consisting of 10 independent segments, each of length 100.

\begin{figure}
\caption{GGS correctly identifies that there are 10 underlying segments in the data (from the
kink in the plots at $K = 9$).}
\begin{subfigure}{0.5\textwidth}
  \centering
  \includegraphics[width=.99\linewidth]{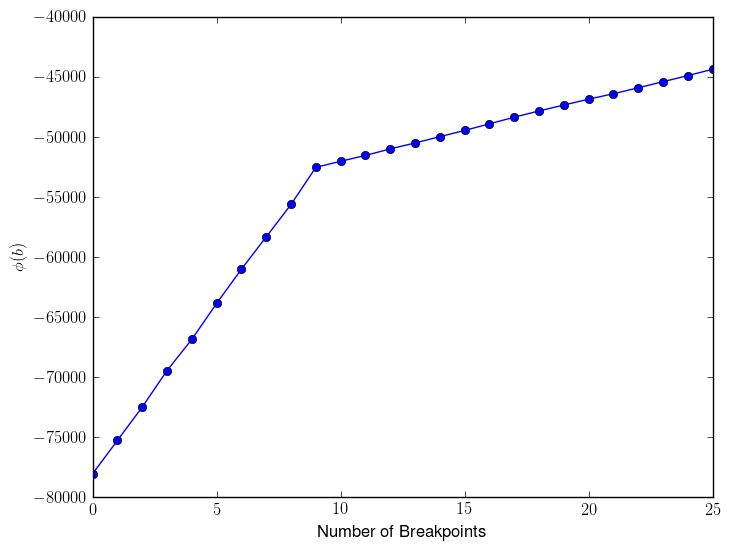}
  \caption{Objective vs.\ breakpoints for $\lambda = 10$.}
  \label{fig:synthFig1}
\end{subfigure}%
\begin{subfigure}{.5\textwidth}
  \centering
  \includegraphics[width=.99\linewidth]{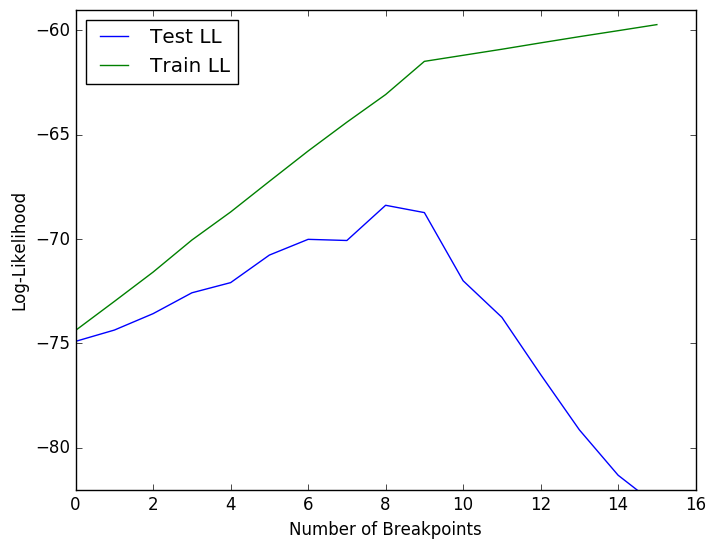}
  \caption{Average train/test log-likelihood for $\lambda = 10$.}
  \label{fig:synthFig2}
\end{subfigure}
\label{fig:synthetic}
\end{figure}

\paragraph{Results.} We run both GGS and the left-to-right HMM 
on this data set. For GGS, we immediately notice a kink in the objective at $K = 9$, 
as shown in Figure \ref{fig:synthFig1},
indicating that the data should be split into $K+1 = 10$ segments.
We use cross-validation to choose an appropriate value of $\lambda$, which yields $\lambda=10$.
We plot the train and test set log-likelihoods at $\lambda=10$ in Figure
\ref{fig:synthFig2}. (Similar to the Wikipedia text example, though, the breakpoint locations are relatively 
robust to our selection of $\lambda$. GGS returns identical breakpoints for any 
$\lambda$ between $10^{-3}$ and $10^3$). 
For this value of $\lambda$ (and thus for the whole range of $\lambda$
between $10^{-3}$ and $10^3$), we split the data perfectly, 
identifying the nine breakpoints at their \emph{exact} locations.

Left-to-right HMMs have various methods for determining the number of segments,
such as AIC or BIC.
Here we instead simply use the correct number of segments, and initialize
the transition matrix as its true value.  Note 
that this is the best-case scenario for the left-to-right HMM. Even 
with this advantage, the left-to-right HMM 
struggles to properly split the time series. Whereas GGS correctly identifies 
$[100, 200, 300, 400, 500, 600, 700, 800, 900]$ as the breakpoints, 
the left-to-right HMM gets at least one breakpoint completely
wrong (and splits the data at, for example, $[100, 200, 300, 400, 500, 600, 663, 700, 800]$).

These results are consistent. In fact, when
this experiment was repeated 100 times (with different randomly generated data), GGS 
identified the correct breakpoints every single time. We also note that GGS is robust to $n$ 
(the dimension of the data), $K$ (the number of breakpoints), and $T/(K+1)$ (the average 
segment length), perfectly splitting the data at the exact breakpoints for all tests across 
at least one order of magnitude in each of these three parameters. On the other hand, 
the 100 left-to-right HMM experiments correctly labeled on average just 7.46 of the nine true breakpoints 
(and never more than eight). 
In these HMM experiments, instead of ten segments of length 100, the shortest segment had an average length
of 26, and the longest segment had an average length of 200.
Additionally, the left-to-right HMM struggles as the 
parameters change, performing even worse when $K$ increases and when $T/(K+1)$ is small compared to $n$ 
(though formal analysis of the robustness of left-to-right HMMs is outside the scope of this paper).
This comes as no surprise, because finding the global maximum among all local maxima of the likelihood function 
for an HMM with many states is known to be difficult problem \cite[Chapter 1.4]{cappe2005inference}.
Therefore, as shown by these these experiments,
GGS appears to outperform left-to-right HMMs in this setting.

\clearpage
\section{Summary}\label{s-summary}

We have analyzed the problem of breaking a
multivariate time series into segments, where the data in each
segment could be modeled as
independent samples from a multivariate Gaussian distribution. 
Our greedy Gaussian
segmentation (GGS) algorithm is able to approximately maximize the 
covariance-regularized log-likelihood in an efficient manner, easily scaling to 
vectors with dimension over 1000 and time series of any length. 
Examples on both small and large data sets yielded useful insights. 
Our implementation,
available at \url{https://github.com/cvxgrp/GGS}, can be used to
solve problems in a variety of applications. For example, the 
regularized parameter estimates obtained by GGS could be used as 
inputs to portfolio optimization, where correlations between 
different assets play an important role when determining optimal holdings. 

\section*{Acknowledgments}
This work was supported by DARPA XDATA, DARPA SIMPLEX, and Innovation Fund Denmark under Grant No.\ 4135-00077B. 
We are indebted to the anonymous reviewers who pointed us to the global solution
method using dynamic programming, and suggested the comparison with
left-to-right HMMs.

%DARPA SIMPLEX Award: N66001-15-C-4042

%As such, there are many benefits of GGS, both in terms of
%driving future research and by directly applying it to real-world 
%problems in various industries.
%
%We define the
%\emph{segmented Gaussian model} (SGM), a covariance-regularized
%maximum likelihood problem, to find the locations of the breakpoints. 
%Choosing these locations
%is a combinatorial optimization problem that is in general difficult to solve,
%so we propose a scalable heuristic that approximately solves it, which we
%call the \emph{greedy Gaussian segmentation} (GGS) algorithm. This
%algorithm is fast, performs well in practice, and yields a locally optimal
%1-OPT solution. We propose a scheme to choose reasonable values
%for the problem's two hyperparameters, and discuss several extensions 
%of this work. Our GGS solver
%is available at (TODO - LINK). Overall, our simple yet efficient algorithm
%can be used for a variety of applications, and has several potential
%extensions to drive future research and to increase the 
%practical benefits of this work.
%The regularized parameter estimates obtained by GGS
%could be useful, for example, as input to portfolio optimization.

\clearpage
\bibliography{ggs}

\end{document}